\documentclass[12pt]{amsart}

\def\@typesizes{%
       \or{5}{6.5}\or{6}{7.5}\or{7}{8.5}\or{8}{11}\or{9}{12}%
       \or{10}{13}% normalsize
       \or{\@xipt}{14}\or{\@xiipt}{15}\or{\@xivpt}{18}%
       \or{\@xviipt}{20}\or{\@xxpt}{24}}

\usepackage{amsthm}
\usepackage{amsmath}
\usepackage{amssymb}
\usepackage{graphicx}
\usepackage{enumerate}
\usepackage{color}
\allowdisplaybreaks
\numberwithin{equation}{section}
\numberwithin{figure}{section}

\theoremstyle{plain}
\newtheorem{theorem}{ Theorem}[section]
\newtheorem{proposition}[theorem]{ Proposition}
\newtheorem{lemma}[theorem]{ Lemma}
\newtheorem{corollary}[theorem]{ Corollary}
\newtheorem{example}[theorem]{ Example}
\newtheorem{remark}[theorem]{ Remark}
\newtheorem{definition}[theorem]{ Definition}
\newtheorem{conjecture}{ Conjecture}

\numberwithin{equation}{section}

\def\BET{\begin{theorem}}
\def\ENT{\end{theorem}}
\def\BEP{\begin{proposition}}
\def\ENP{\end{proposition}}
\def\BEL{\begin{lemma}}
\def\ENL{\end{lemma}}
\def\BEC{\begin{corollary}}
\def\ENC{\end{corollary}}
\def\BEE{\begin{example} \rm}
\def\ENE{\end{example}}
\def\BER{\begin{remark} \rm}
\def\ENR{\end{remark}}
\def\BED{\begin{definition} \rm}
\def\END{\end{definition}}
\def\BECJ{\begin{conjecture}}
\def\ENCJ{\end{conjecture}}

\def\bea{\begin{eqnarray}}
\def\eea{\end{eqnarray}}

\def\beas{\begin{eqnarray*}}
\def\eeas{\end{eqnarray*}}

\def\beq{\begin{equation}}
\def\eeq{\end{equation}}

\def\beal{\begin{align*}}

\def\eeal{ \end{align*} }

\def\roweq{\nonumber \\ &=& }
\def\rowleq{\nonumber \\  & \leq & }

\def\rowpl{\nonumber \\ &  \ \ + & }

\def\bfQ{{\bf Q}}

\def\bbC{{\mathbb C}}

\def\bbN{{\mathbb N}}

\def\bbR{{\mathbb R}}

\def\bbZ{{\mathbb Z}}

\def\cA{{\mathcal A}}

\def\cH{{\mathcal H}}
\def\cI{{\mathcal I}}

\def\cL{{\mathcal L}}
\def\cM{{\mathcal M}}

\def\cP{{\mathcal P}}

\def\cR{{\mathcal R}}

\def\cT{{\mathcal T}}

\def\cX{{\mathcal X}}

\def\sF{{\sf F}}

\def\wP{{\widetilde P_\mu}}

\def\sft{\mathsf{t}}
\def\tr{\mathsf{tr}}

\def\ef{\eqref}

\begin{document}

\title[Elliptic functions, Floquet transform, Bergman spaces]{Elliptic functions, Floquet transform and Bergman spaces on doubly periodic domains}

\author{Jari Taskinen$^1$}

\author{Zhan Zhang$^2$}

\thanks{The research of the first named author was partially supported by the Academy of Finland
grant no. 359563}

\begin{abstract}
We study Bergman spaces $A^2(\Omega)$, their kernels and  Toeplitz operators on 
unbounded, doubly periodic domains $\Omega$ in the complex plane.
We establish the mapping  properties of the Floquet transform  operator defined in $A^2(\Omega)$
and derive a general formula connecting the Bergman kernel and projection of the domain 
$\Omega$ to a kernel and projection on the bounded periodic cell $\varpi$. As an application, we 
prove, for Toeplitz  operators $T_a$ with  doubly periodic symbols, a spectral band formula, which describes 
the spectrum and essential spectrum of $T_a$ in terms of the spectra of a family of Toeplitz-type operators 
on the cell $\varpi$.

Technical challenges arise from the fact that double quasiperiodic boundary conditions have to be
taken into account in the definitions of the spaces and operators on the periodic cell $\varpi$. 
This requires novel operator theoretic tools, which are based on modifications of certain elliptic functions, e.g. the Weierstrass $\wp$-function. 
\end{abstract}

\maketitle

{\small
\noindent Addresses:

\smallskip

\noindent 
$^1$ Department of Mathematics and Statistics,
 University of Helsinki, P.O. Box 68, 00014 Helsinki, Finland

\noindent {\it email} jari.taskinen@helsinki.fi

\smallskip

\noindent $^2$ Department of Mathematics and Statistics,
 University of Helsinki, P.O. Box 68, 00014 Helsinki, Finland

\noindent {\it email} zhan.zhang@helsinki.fi}

\section{Introduction.}
\label{sec1}

In this paper we consider Bergman space $A^2(\Omega)$ and the Bergman 
projection $P_\Omega$ on planar domains $\Omega \subsetneq \bbC$, which have the special 
geometry of being  periodic in  two directions. 
% the domain is obtained as 
%the union of infinitely  many translated copies of the bounded periodic 
%cell $\varpi \subset \bbC$. 
For simplicity, we require that
\bea
z \in \Omega \ \Rightarrow \ \ z + m \in\Omega \ \mbox{for all} \ 
m  \in \Lambda %%\bbL_2 
= \{ z \in \bbC \, : \, {\rm Re}\,z \in \bbZ , \ {\rm Im}\, z \in \bbZ \} .
\label{1.1}
\eea 
%see below for more details of the definitions. 
Following the paper \cite{T1}, where the singly periodic case was treated, the 
aim is to apply Floquet transform 
techniques to describe the connection of the Bergman   projection and kernel on $\Omega$ with projections and kernels on
the periodic cell $\varpi$. We will apply these considerations to study the spectra of Toeplitz operators 
with doubly periodic symbols. 

Given a domain $\Omega$ in the complex plane $\bbC$, we denote by $L^2(\Omega)$ the usual Lebesgue-Hilbert space with respect to 
the (real) area measure $dA$ and by  $A^2(\Omega)$ the corresponding Bergman space, which
is the closed subspace consisting of analytic functions.   The Floquet transform
is defined  for $f \in A^2(\Omega)$ by
\bea 
\sF f(z,\eta)&=& \widehat f (z ,\eta ) =  \frac{1}{2\pi}\sum_{m\in \Lambda}
e^{-i \eta_1m_1 - i \eta_2 m_2  } f (z + m ), \label{2.2}
\eea
where $ z \in \varpi$ and $\eta = (\eta_1, \eta_2) \in \bfQ = [-\pi,\pi]^2$ is the so called Floquet-parameter or
quasimomentum; here and later, it will be convenient to denote $m_1 = {\rm Re}\, m$, $m_2 = \ {\rm Im}\, m$ for a number $m \in \Lambda$.
By using the canonical identification of $\bbC \cong \bbR^2$, this definition 
coincides with the standard definition of the Floquet transform in periodic subdomains 
of $\bbR^2$. In particular, the well-known properties of Fourier series (see also \cite{Ku}, Section 4.2) 
imply that $\sF$ is a  continuous (isometric) mapping from $ L^2(\Omega)$ onto $L^2(\bfQ ; L^2(\varpi))$,
which is the vector valued $L^2$ space, or the space of $L^2$-Bochner integrable functions $g :\bfQ \to 
L^2(\varpi)$; see \cite{Hyt}.  If  $g \in L^2(\bfQ; L^2(\varpi))$, we also denote
\bea
\sF^{-1} g (z ) & = &  \frac1{2 \pi} \int\limits_{\bfQ}
e^{ i [{\rm Re} z] \eta_1 + i [{\rm Im} z] \eta_2} g(z -  [{\rm Re}\,z ] - 
i [{\rm Im} z] , \eta ) d \eta ,  \ \ z \in \Omega  .
%%e^{ i ( [{\rm Re}\,z] - z ) \eta} ?? 
\label{2.4}  
\eea 

We refer to \cite{Kbook}, \cite{Ku} for an introduction of the Floquet transform in the study of the Schr\"odinger equation with
periodic potentials and to \cite{Na1}, \cite{NaPl} for its use in periodic elliptic spectral problems. See also \cite{T1}
for a slightly more thorough introduction to the topic in the setting of analytic function spaces.  
Taking the restriction of $\sF$ to the subspace of analytic functions, we conclude that $\sF$ is a unitary  mapping 
from $ A^2(\Omega)$ into $L^2(\bfQ ;  L^2(\varpi))$, but in order to achieve the bijectivity, there remains
to determine the image $\sF(A^2(\Omega))$ in the Bochner space.  In \cite{T1},  the first named author established the
mapping properties of $\sF$ and solved this problem in the case of complex domains $\Pi$, which are periodic, say, in the
direction of the real axis and bounded in the imaginary direction and thus contained in a strip. 

The aim of this paper is to extend many of the results of  \cite{T1}, \cite{T2} to the doubly periodic case, but having a 
closer look at the arguments of these papers, one observes that two essential technical tools of the citations fail in
the present setting. To characterize the range of the Floquet transform in the Bergman space case on singly periodic domains, \cite{T1},
it was necessary 
to approximate a Bergman function $f \in A^2(\Pi)$ by a function $\varphi_\varepsilon f$, where $\varphi_\varepsilon (z) = e^{-\varepsilon z^2}$, $z \in \Pi$, $\varepsilon>0$, is a rapidly decreasing function as $|z| \to \infty$ in $\Pi$. Indeed,
the function $\varphi_\varepsilon$ has Gaussian decay at infinity, if $z$ belongs to a strip parallel to the real axis,
but loses this property, if $z$ is allowed to belong to a doubly periodic domain, which is automatically unbounded in
the imaginary direction. The same problem appears when trying to generalize the results of \cite{T2}, where the operator
$J_{\nu,\mu} : f \mapsto e^{i (\nu - \mu)z} f$ was used to switch the Floquet parameters $\nu,\mu \in [-\pi,\pi]$
in the quasiperiodic boundary conditions.  

In Section \ref{sec2} we construct, by using the theory of elliptic (doubly periodic meromorphic) functions, the tools
that allow us to adapt the methods of \cite{T1}, \cite{T2} to the doubly periodic case. These consist of the functions 
$\psi^{(\rho)}$ and $\psi_\eta$, see  Lemmas \ref{lem1.3} and \ref{lem1.6}, which replace the functions 
$e^{-\varepsilon z^2}$ and $ e^{i \eta z}$ used in the singly periodic case. 

The remaining sections are devoted to the formulation of the main results and their  proofs. In Section \ref{sec3}, Theorem \ref{th2.3}, 
we characterize the image $\sF(A^2(\Omega))$ in $L^2(\bfQ;L^2(\varpi))$ and in  Theorem \ref{th3.3} of Section \ref{sec4} we present the relation of the Bergman kernels in the
domains $\Omega$ and $\varpi$. In  Section \ref{sec5} we apply the results of the preceding sections to study
spectra of Toeplitz operators $T_a$ with periodic symbols $a$  on the space $A^2(\Omega)$.  Theorem \ref{th8.3} contains a proof for the spectral band 
formula, which characterizes the spectrum and essential spectrum of $T_a$ in terms of a family of spectra of $\eta$-dependent Toeplitz-type
operators in the periodic cell. 

We denote by $P_\Omega$ the orthogonal projection from $L^2(\Omega)$ onto $A^2(\Omega)$. It can 
always written with the help of the Bergman kernel $K_\Omega : \Omega \times \Omega \to \bbC$,
\bea
P_\Omega f(z) = \int\limits_\Omega K_\Omega (z,w) f(w) dA(w)
\eea
and the kernel has the properties that $K_\Omega (z, \cdot) \in L^2(\Omega)$ 
for all $z$ and $K(z,w) = \overline{ K(w,z)}$ for all $z,w \in \Omega$. 
See e.g. \cite{K1} for a proof of these assertions. 

The basic assumptions on the domain $\Omega$ and its periodic cells are as follows. However, we will pose one more assumption, related
to the theory of elliptic functions, see $(\cA)$ in the next section. Here and later, cl$(A)$ denotes the closure of the set $A$. 

\BED  \label{def1}
We denote  by  $\varpi$ the periodic cell, which is a domain in $\bbC$ such that 

\smallskip \noindent
$(i)$ %%\bea
$\varpi \subset
%B\big( \bfc , %(1+i)/2,  1- 1/\sqrt{2}\big) \subset 
Q= (0,1) \times (0,1) \subset \bbR^2 \cong \bbC, %%\label{1.0}
$ %% \eea 

\smallskip \noindent
$(ii)$ the  boundary  $\partial \varpi$ contains the  boundary $\partial Q$, and

\smallskip \noindent
$(iii)$ the complement $Q \setminus {\rm cl} (\varpi) $ contains an open set

\smallskip
 
We denote the translates of $\varpi$ by $\varpi_m = \varpi + m$ for all $m \in \Lambda$.
%%(Recall that we  denote $m_1 = {\rm Re}\, m$, $m_2 = \ {\rm Im}\, m$ for all  $m \in \Lambda$.)
The periodic domain $\Omega$ is defined  as the interior of the set
\beas
\bigcup_{m \in \Lambda} {\rm cl}( \varpi_m ).   %%\label{1.2}
\eeas
\END

We emphasize that contrary to \cite{T1}, we do not need to add assumptions on the smoothness of the boundary $\partial \Omega$.
Anyway, the periodic domain $\Omega$ is always infinitely connected, due to assumption $(iii)$ of Definition \ref{def1}.

\bigskip
The following notation will be used throughout the paper. We write $C$, $C'$, $\ldots$,
(respectively, $C_\eta$, $C'_\eta$, $\ldots$ etc.) for positive constants independent  
of the functions or variables in the given inequalities (resp. depending only on a 
parameter  $\eta$ etc.), the values of which may vary from place 
to place. For $x \in \bbR$, $[x]$ denotes the largest integer not larger than
$x$. We write $\bbZ= \{ 0, \pm1, \pm2, \ldots \}$. 
Subsets of the complex plane will often be described using the real planar coordinates,  
%canonical identification $\bbC \cong \bbR^2$, 
for example, rectangles in $\bbC$ are described as the sets $[a,b]\times [c,d]$. 
If $x$ is a point in $ \bbC$ or $\bbR^n$ and $r >0$, then $B(x,r)$ 
denotes the Euclidean ball with center $x$ and radius $r>0$. Moreover,
cl$(A)$ denotes the closure of a set $A$. 

Given a domain $D \subset \bbC $, we denote by $\Vert f \Vert_D$ and $( \cdot | \cdot)_D$ the norm 
and inner product of $L^2(D)$, which is the $L^2$ space with respect to the (real) area measure.
In general, the norm of a Banach space $X$ is denoted by $\Vert \cdot \Vert_X$. 
Given an interval $I \subset \bbR$ and a Banach space $X$ we denote  by $L^2(I;X)$ the 
space of vector valued, Bochner-$L^2$-integrable functions on $I$ with values  in  $X$, 
endowed  with the norm 
\beas
\Vert f \Vert_{L^2(I;X)}=\Big(\int\limits_I 
\|f( t )\|_X d t \Big)^{1/2}\, .
%\label{1.7}
\eeas 
If $X = L^2(D)$ for some domain $D$, then $L^2(I;X)$ is a Hilbert space 
endowed with the  inner product $\int_I (f (t)| g (t))_D dt  $.
See \cite{Hyt} for the theory of Bochner spaces. 

Given a Banach space $X$, $\cL(X)$ stands for the Banach space of bounded
linear operators $X \to X$. The operator norm of $T \in \cL(X)$ is  denoted just by
$\Vert T \Vert$ or by $\Vert T \Vert_{X \to X}$, if it is necessary to specify the 
domain or target spaces. The identity operator on $X$ is written as $I_X$ or by
$I_D$, if $X=L^2(D)$.  
For an operator  $T \in \cL(H)$, where $H$ is a Hilbert space,  $\sigma (T)$, $\sigma_{\rm ess} (T)$ and 
$\varrho(T)$ stand for the spectrum, essential spectrum and resolvent set of the 
operator $T$. The resolvent (operator) of  $T \in \cL(H)$ is denoted by 
$R_\lambda(T) = (T -\lambda I_H)^{-1}$, where $\lambda \in \varrho(T)$.

\section{On elliptic functions in the periodic domain  $\Omega$.}
\label{sec2}

Our aim in this section is to prove Lemmas \ref{lem1.3} and \ref{lem1.6} which provide the technical tools needed for
the proofs of the main results in later sections. The tools consist of modifications of certain elliptic functions. 
Recall that an elliptic function  $\varphi$ is by definition a meromorphic and doubly periodic function
in the plane $\bbC$, see for example  \cite{FB}, Chapter 5, or \cite{NP}, Chapter 14, and others. It ought to be recalled that, except for the constant functions, there does not exist doubly
periodic functions which are analytic in the entire complex plane (Theorem 1, Chapter 14 of \cite{NP}), and the use of meromorphic ones is thus necessary for our
purposes.   Here, in view of our choice of the domain $\Omega$, we will assume the periodicity
\bea
\varphi(z) = \varphi (z+1) = \varphi (z +i) .  \label{1.4}
\eea
According to the classical theory, any elliptic function satisfying \ef{1.4} has a finite number of zeros in the cell $[0,1) \times [0,1) 
\subset \bbC$ and equally many poles there as well, when counted by taking into account the zero and pole orders (Theorem 5, 
Chapter 14 of \cite{NP})

Let us make the following assumption on the domain $\Omega$.

{\it $(\cA)$ There exists an elliptic function $\varphi$ with periodicity \ef{1.4}, such that the zeros of $\varphi$ occur at points $z \in \bbC \setminus {\rm cl}( \Omega)$ and such that $\varphi$ has a representation
\bea
\varphi(z) = \sum_{m \in \Lambda} \phi(z - m)   \label{1.6}
\eea 
for some function $\phi$, which is meromorphic in $\bbC$ and has only finitely many poles, all situated in $Q$ 
%%$\varpi$, IN Q I THINK!?!?!?!
and which satisfies
\bea
\sup_{z \in \Omega \setminus {\rm cl}( Q) }(1 +|z|)^b |\phi(z)| =: d_\phi <  \infty   \label{1.7}
\eea
for some constant $b> 2$.
}

We denote by $S \subset Q$ the finite set of the poles of the function $\phi$. 
Note that \ef{1.7} implies the absolute convergence of the series \ef{1.6}, and there actually holds a stronger statement, see 
Lemma \ref{lem1.2}.

\BEE \label{ex2.1}
We show that $(\cA)$ holds, if $Q \setminus {\rm cl}(\varpi)$ contains the three points $\alpha + \frac12$, $\alpha + \frac{i}2$
and $\alpha + \frac{1+i}2$ for some $\alpha$ with Re\,$\alpha \in (0,\frac12)$ and Im\,$\alpha \in (0,\frac12)$. Indeed, 
the Weierstrass $\wp$-function is defined by
\bea
\wp(z) = \frac{1}{z^2 }+ \sum_{m \in \Lambda \setminus \{0\} } \Big( \frac{1}{(z-m)^2} - \frac{1}{m^2 } \Big) ,
\label{1.8}
\eea 
and its well-known properties include the following (\cite{NP}, Chapter 14) : 

\smallskip \noindent $(i)$ $\wp$ is a meromorphic function in $\bbC$ with poles of order 2 at the points $m\in \Lambda$,

\smallskip \noindent $(ii)$ the zeros of $\wp$ occur at the points $\frac12(1+i) + m$, $m \in \Lambda$, %%[[Duke-Imamoglu p. 902]]
and

\smallskip  \noindent $(iii)$ $\wp$ is doubly periodic 
with periodicity \ef{1.4}. %% more precisely, $\wp(z) = \wp(z+1) = \wp (z+i)$ for
%%all $z \in \bbC \setminus \Lambda$.

\smallskip

\noindent The representation \ef{1.8} is not of the form \ef{1.6}, but the derivative $\wp'$ can be calculated by termwise differentiation, which yields
\beas
\wp'(z) = \sum_{m \in \Lambda}  \frac{-2}{(z-m)^3} ,  %\label{1.10}
\eeas
i.e., a function of the form \ef{1.6}. It is known that the zeros of $\wp'$ occur at the points 
$\frac12$, $\frac{i}{2}$ and $\frac12(1+i)$ and their translates, see the references above, or \cite{DI}, Section 1. %%WIKIPEDIA WEIERS. P-FUNCTION. 
We define $\phi(z) = -2(z- \alpha)^{-3}$ and $\varphi(z) = \wp'(z - \alpha)$, where $\alpha$ is a above. It is plain that with these
functions, the domain $\Omega$ satisfies condition $(\cA)$. 
\ENE

From now on we assume that  $\Omega$ is a doubly periodic domain  satisfying condition $(\cA)$ with functions
$\varphi$ and $\phi$. In view of Example \ref{ex2.1}, this is only a minor restriction of generality. 
%%[[MAYBE SOME WORDS ABOUT EEROS CONSTRUCTION, WHICH DOES NOT SUFFICE.]]

\BER  \label{rem1.2}
The function $\varphi$ is bounded from below on $\Omega$: there exists a constant $C> 0$ such that
$\inf_{z \in \Omega} |\varphi(z)| \geq C$.  Namely, by the  basic assumptions, the zeros of $\varphi$ are outside the set 
$ {\rm cl} ( \Omega)$. The claim follows from the periodicity of $\varphi$, which allows us to restrict the proof of the claim
into a compact subset of $\bbC$.  
\ENR

We will need two modifications of the function $\varphi$. The first one is given in Lemma \ref{lem1.3} and it gives a replacement 
$\psi^{(\rho)}$ of the function $\varphi_\varepsilon$
of the singly periodic case (see Section \ref{sec1}).  We first consider the following statement.

\BEL  \label{lem1.2}
If $K$ is a compact set contained in $\bbC \setminus \bigcup_{m \in \Lambda} (S+m)$, then, for every $\delta \in (0,b-2)$, 
\bea
\sum_{m \in \Lambda} |m|^{b-2- \delta} \sup_{z \in K \cap \Omega} |\phi(z-m)| < \infty. \label{1.11}
\eea
\ENL

Proof. %%If $m \in \Lambda$ and $z \in \varpi_m$, then $|z|$ is proportional to $|m|$. Moreover, 
Given a compact $K$   as in the assumption, 
there are constants $C, C' > 0$ such that $|m| \leq C |z- m|$ for all $z \in K$, if  $|m| \geq C'$.
We obtain 
\bea
& & \sum_{\stackrel{\scriptstyle m \in \Lambda}{|m| \geq C'}} |m|^{b-2- \delta} \sup_{z \in K \cap \Omega} |\phi(z-m)|
= \sum_{\stackrel{\scriptstyle m \in \Lambda}{|m| \geq C'}} |m|^{-2- \delta} \sup_{z \in K \cap \Omega} |m|^b |\phi(z-m)|
\rowleq
C\sum_{\stackrel{\scriptstyle m \in \Lambda}{|m| \geq C'}} |m|^{-2- \delta} \sup_{z \in K \cap \Omega} |z- m |^b |\phi(z-m)|
\label{1.11a}
\eea
By  \ef{1.7}, we have here  
\bea
|z- m |^b |\phi(z-m)| \leq  d_\phi   \label{NEW1}
\eea
for all $z \in K$ except possibly for $z $ belonging to the cell $\varpi_m$, since $\phi$ has poles in $Q$. However, if $z \in \varpi_m
\cap K$, then $|z-m| \leq 1$ and $|\phi(z-m)|$ is still bounded by a constant depending on $K$ only, since we are assuming that the 
distance of $K$ from the set $S+ m$ is positive. Hence, \ef{NEW1} holds for all $K$, with possibly a 
larger constant on the right hand side.  We conclude that \ef{1.11a} is bounded by
\beas
C \sum_{m \in \Lambda} |m|^{- 2- \delta} < \infty. 
\eeas
The remaining finitely many terms of \ef{1.11} with $|m| < C'$ are uniformly bounded, by the assumption 
of the lemma. \ \ $\Box$

\bigskip

The construction and important properties of $\psi^{(\rho)}$ are contained in the next lemma. It uses 
the elliptic function $ \varphi$  associated with the domain $\Omega$. Recall that $b>2$ was determined 
in \ef{1.7}.

\BEL  \label{lem1.3}
Let $\varphi$, $\phi$ be as in \ef{1.6}, let $\rho  \in (0,1]$ and define  the functions 
\bea
\varphi^{(\rho)} (z) = \sum_{m \in \Lambda} 2^{ -\rho |m|  } \phi(z - m), \ z \in \Omega,   \ \ \ 
\mbox{and} \ \ \ \psi^{(\rho)}  = \frac{\varphi^{(\rho)} }{\varphi}.   \label{1.12}
\eea

\smallskip \noindent
$(i)$ For every $\rho \in (0,1]$, the functions $\varphi^{(\rho)}$ and $\psi^{(\rho)}$ 
are meromorphic in $\bbC$, and $\psi^{(\rho)}$ is analytic in $\Omega$. 

\smallskip \noindent
$(ii)$ There exists a constant $C >0$ such that 
\bea
\sup_{z \in \Omega, \rho \in(0,1] } |\psi^{(\rho)}(z)| \leq C . \label{1.12h}
\eea
Moreover, for every $M \geq 1$, there holds %%has the property that 
$\psi^{(\rho)} \to 1 $ uniformly on $\big( [-M, M]\times [-M, M] \big) \cap \Omega$ as $\rho \to 0$. 

\smallskip \noindent
$(iii)$
Given $f \in L^2 (\Omega)$, there also holds
\bea
\Vert f \psi^{(\rho)} - f  \Vert_{\Omega} \to 0 \ \ \ \mbox{as} \ \rho \to 0 .  \label{1.14}
\eea

\smallskip \noindent
$(iv)$ For all compact $K \subset \bbC$ and $m \in \Lambda$, $\rho \in (0,1]$, we have 
\bea
\sup_{z \in K \cap \Omega} |\psi^{(\rho)} (z+m)| \leq \frac{C_K}{\rho^b(1 + |m|)^b}. 
%\ \ \ \ \mbox{for all} \ .
\label{1.15}
\eea
\ENL 

Proof. $(i)$ Let us fix $\rho \in (0,1]$. First, due to Lemma \ref{lem1.2}, the series in \ef{1.12} converges
absolutely and uniformly on every compact set contained in  $\bbC \setminus \bigcup_{m \in \Lambda} (S+m)$, hence, 
the function $\varphi^{(\rho)}$ is analytic in $\bbC \setminus \bigcup_{m \in \Lambda} (S+m)$. Since $\phi$ is
meromorphic and the poles are contained in $Q$, it follows that $\varphi^{(\rho)}$ is meromorphic in $\bbC$. Also, since $\varphi$ is meromorphic, 
$\psi^{(\rho)}$ is meromorphic in $\bbC$, too.

Let us show  the analyticity of $\psi^{(\rho)}$ in $\Omega$. Since $\varphi$ does not have zeros in a neighborhood of ${\rm cl}( \Omega) $, 
the function $\psi^{(\rho)}$
is analytic everywhere in $\Omega$ except possibly at the poles of $\varphi^{(\rho)}$. So, let 
$z_0 \in \Omega$ be a pole  of $\varphi^{(\rho)}$ of order $N \in \bbN$,   belonging to cl$(\varpi_{m_0})$ for some $m_0 \in \Lambda$.
We have 
\bea
\varphi^{(\rho)}(z) = \sum_{n=1}^N \frac{a_n^{(\rho)}}{(z-z_0)^n} + \alpha^{(\rho)}(z)   \label{1.16}
\eea
for some complex numbers $a_n^{(\rho)}$ and some function $\alpha^{(\rho)}$ which is analytic in a neighborhood $U$ of $z_0$.

If $z_0 \in \partial ( Q + m_0))$, we can choose $U$ so small that $|z-m| \geq 1/2 $ for all $z \in U$ 
with  $|m - m_0| \geq 2$. By \ef{1.7}, $|\phi(z-m)| \leq C_\phi (|m-m_0|+1)^{-b} $ for 
all such $m$ and all  $z \in U$, hence, 
\bea
& & \sum_{|m - m_0| \geq 2} 2^{-  \rho|m|} |\phi(z-m)|   %%\sum_{|m - m_0| \geq 2}  |\phi(z-m)|   
\leq \sum_{|m - m_0| \geq 2} C_\phi (|m - m_0|+1)^{-b}   \leq  C_\phi'  \label{1.20}
\eea
for $z \in U$. We conclude that the pole  of  $\varphi^{(\rho)}$ at $z_0$ is determined only by the 
terms with   the functions 
$\phi(z-m)$ with $|m - m_0| \leq  1$ in the series \ef{1.12}. But this is not possible, since the poles 
of $\phi$ were assumed to be in (the interior of) $Q$, not on the boundary, see condition $(\cA)$. Hence,
$z_0 \notin \partial ( Q + m_0)$.

We can thus choose $U$ such that  dist\,$(U, \partial ( Q + m_0)) \geq \delta $ for some $\delta > 0$,
and we obtain that $|z-m| \geq \delta$ for all $z \in U$, $m \not=m_0$. Again, 
$|\phi(z-m)| \leq C (|m-m_0|+1)^{-b} $ for some constant depending on $\phi$ only and for all $z \in U$ 
and $m \not=m_0$, which yields  
\bea
& & \sum_{m \not= m_0} 2^{-  \rho|m|} |\phi(z-m)|   \leq
\sum_{m \not= m_0}  |\phi(z-m)|   \leq \sum_{m \not= m_0} C (|m - m_0|+1)^{-b}   \leq  C_\phi'  \label{1.20}
\eea
for some positive constant $C_\phi'$, independent of $\rho$ (we will need this later). The pole  of 
$\varphi^{(\rho)}$ at $z_0$ is thus determined only by the term with the function $\phi(z-m_0)$ in the 
series \ef{1.12}. More precisely, we deduce by comparing \ef{1.12}, \ef{1.16} and \ef{1.20} that
\bea
\phi(z- m_0) = \sum_{n=1}^N \frac{ 2^{\rho |m_0|} a_n^{(\rho)}}{(z-z_0)^n} + \beta^{(\rho)}(z)  \label{1.22}
\eea 
for some function $\beta^{(\rho)}$ which is analytic  in a neighborhood of $z_0$. 

Due to \ef{1.7} and an estimation similar to \ef{1.11a}, the bound \ef{1.20} holds also in the case $\rho =0$, which, in view of \ef{1.6} and \ef{1.22} implies
that also 
\bea
\varphi(z) =  \sum_{n=1}^N \frac{2^{\rho |m_0|} a_n^{(\rho)}}{(z-z_0)^n} + \gamma^{(\rho)}(z)  \label{1.23}
\eea 
for another function $\gamma^{(\rho)}$ analytic in a neighborhood of $U$. 

We conclude that in the definition of $\psi^{(\rho)}$, the poles of $\varphi^{(\rho)}$
and $\varphi$ at the point $z_0$ cancel each other, i.e., the function $\psi^{(\rho)}$ is analytic in $U$. 
Since $z_0$ was an arbitrary pole of $\psi^{(\rho)}$, the function  $\psi^{(\rho)}$ must be analytic in 
$\Omega$. 

$(ii)$ We first claim that,  for any $m \in \Lambda$,
\bea
z \mapsto \frac{\phi(z -m )}{\varphi(z)}  \label{1.24}
\eea
is a bounded function in $\Omega$. Namely, by formulas \ef{1.22} and \ef{1.23}, the function 
$ \phi(z-m)/\varphi(z)$ is analytic in a neighborhood of 
the finitely many poles of $\phi( z - m)$. Outside this neighborhood, the function $\phi( z - m)$ is uniformly bounded,
due to the assumption \ef{1.7}. Also, by Remark \ref{rem1.2}, $|\varphi(z)| \geq C$ for a positive constant $C$, for all 
$z \in \Omega$. Putting these observations together proves the claim.  

We next show that \ef{1.12h} holds. 
% We first note that the functions $|\psi^{(\rho)}|$ have a uniform upper bound in $\Omega$, i.e.,
%there holds
%\bea
%M_0 = \sup_{ z \in \Omega, \rho > 0} \{|\psi^{(\rho)}(z)| \} <\infty. \label{1.28}
%\eea 
Namely, if $m_0$ is such that $z \in {\rm cl} ( \varpi_{m_0} )$, then we have 
\bea
& & |\psi^{(\rho)}(z) | \leq
\frac{| \phi(z - m_0)|}{|\varphi(z)|} 
+ \frac{1}{|\varphi(z)|} \sum_{\stackrel{\scriptstyle m \in \Lambda}{m \not= m_0}} |\phi(z-m) | ,    \label{1.28u}
\eea
and here, the first term has an upper bound independent of $\rho$ or $z$ by the remark above.
%% in the  beginning of the proof $(ii)$. 
Also, the second term is bounded a constant independent of $\rho$ or $z$,
by Remark \ref{rem1.2}, \ef{1.7} and an argument similar to \ef{1.20}. We obtain \ef{1.12h}.

To prove the uniform convergence, let $M \in \bbN$ and  $\varepsilon>0$ be given. Denoting 
$K = \big( [-M,M] \times [-M, M] \big)  \cap \Omega$, we  use \ef{1.7} (see also the calculation in \ef{1.20}) 
to  choose $N \in \bbN$ such that 
\bea
& & \sup_{z \in K } \sum_{\stackrel{\scriptstyle m \in \Lambda}{|m| > N}} |\phi(z-m)|
\leq %%\rowleq
\sum_{\stackrel{\scriptstyle m \in \Lambda}{|m| > N}}  \sup_{z \in K } |\phi(z-m)|
< \varepsilon.  \label{1.26}
\eea
Then, by \ef{1.6}, 
\beas
& & | \psi^{(\rho)} -  1 | = \frac{1}{|\varphi(z)|} \Big| \sum_{m \in \Lambda} 2^{-\rho |m|} \phi(z-m) -
\sum_{m \in \Lambda} \phi(z-m) \Big|
\rowleq    \sum_{\stackrel{\scriptstyle m \in \Lambda}{|m| \leq N}} \big( 1 - 2^{-\rho |m|} 
\big) \frac{|\phi(z-m) |}{|\varphi(z)|} 
\rowpl
\frac{1}{|\varphi(z)|} \Big| \sum_{\stackrel{\scriptstyle m \in \Lambda}{|m| > N}} 2^{-\rho |m|} \phi(z-m)\Big|
+ \frac{1}{|\varphi(z)|} \Big| \sum_{\stackrel{\scriptstyle m \in \Lambda}{|m| > N}} \phi(z-m) \Big| .
\eeas
Here, since $|\varphi(z)| \geq C>0$ for $z \in \Omega$,  the last two terms are bounded by $\varepsilon/C$,  by \ef{1.26}.
Also, the first term can be estimated using the boundedness of the finitely many expressions \ef{1.24},
\beas
\sum_{\stackrel{\scriptstyle m \in \Lambda}{|m| \leq N} } \big( 1 - 2^{-\rho |m|} \big) 
\frac{|\phi(z-m) |}{|\varphi(z)|} \leq C_N \big( 1 - 2^{-\rho N } \big),
\eeas
which is bounded by $\varepsilon$, if $\rho$ is small enough. Hence, $\psi^{(\rho)} \to
1$ uniformly on $\big( [-M,M] \times [-M, M] \big)  \cap \Omega$.

$(iii)$
If $f \in L^2(\Omega)$ and $\varepsilon>0$ is given,  we pick up %%use \ef{1.7} to find 
$M \in \bbN$ so large that
\bea
\int\limits_{\Omega \cap \{ |z| \geq M \}} \!\!\!\!\!\!\!\! (1+  C)^2 |f|^2 dA < \varepsilon,
\label{1.30}
\eea
where $C$ is the constant in \ef{1.12h}. Then, by \ef{1.30},
\bea
& & \int\limits_{\Omega} |f \psi^{(\rho)} - f|^2 dA  
\leq \!\!\!\!\!\!   \int\limits_{\Omega \cap \{ |z| \leq M \} } \!\!\!\!\!\!
%%\stackrel{\scriptstyle x \in \Omega}{|x| \leq M}} 
|1-\psi^{(\rho)}| \,| f|^2 dA
+ \!\!\!\!\!\! \int\limits_{\Omega \cap \{ |z| \geq M \} }  \!\!\!\!\!\!
%%\stackrel{\scriptstyle x \in \Omega}{|x| \geq M}} 
(1+C)^2 |f|^2 dA 
\rowleq
 \Vert f \Vert_{\Omega}^2 \sup_{|z| \leq M } |1-\psi^{(\rho)}(z)|  + \varepsilon .
\eea
The first term on the right hand side can be made smaller than $\varepsilon$, if $\rho$ is small enough, 
by using the already proven claim $(ii)$ of the theorem, which proves the claim $(iii)$. 

$(iv)$ Let $K$ and $m$ be given and $z \in K$. We have
\beas
& & |\psi^{(\rho)}(z+ m)| \leq 
%%\sum_{n \in \Lambda} 2^{-\rho |n|} | \phi(z-n)| \rowleq  
\!\!\!\!\!\!
\sum_{\stackrel{\scriptstyle n \in \Lambda}{|n - m| \leq |m|/2}}\!\!\!\!\!\!  2^{-\rho |n|} 
\frac{|\phi(z+m -n)|}{|\varphi(z)|} 
+ \!\!\!\!\!\! \sum_{\stackrel{\scriptstyle n \in \Lambda}{|n- m| > |m|/2}} \!\!\!\!\!\! 2^{-\rho |n|} 
\frac{|\phi(z+m -n)|}{|\varphi(z)|} 
\nonumber \\
&& =:R_1 + R_2.
\eeas
To estimate $R_1$, the summation index satisfies $|n| \geq |m|/2$, hence, $2^{-\rho |n|} \leq 2^{- \rho |m|/2}$ and
this yields 
\beas
R_1 \leq 2^{- \rho |m|/2} \sum_{n \in \Lambda } \frac{|\phi(z+m-n)|}{|\varphi(z)|} . 
%%\leq C2^{- \rho |m|/2}   \leq \frac{C'}{\rho^b(1+ |m|)^b},
\eeas
Here, the sum is bounded by a constant independent of $z,m$, by an argument similar to that following \ef{1.28u}. Hence,
\beas
R_1 \leq  C2^{- \rho |m|/2} 
\leq \frac{C'}{\rho^b(1+ |m|)^b}. 
\eeas
As for $R_2$,
\bea
& & R_2 \leq \!\!\!\!\!\! \sum_{\stackrel{\scriptstyle n \in \Lambda}{|n- m| > |m|/2}} \!\!\!\!\!\! 2^{-\rho |n|} 
|m-n|^{-b} \frac{|m-n|^{b} |\phi(z+m-n)|}{|\varphi(z)|} 
\rowleq |m/2|^{-b}  \sup_{w \in K, n\in \Lambda} \frac{|n|^{b} |\phi(w+n)|}{|\varphi(z)|}
\sum_{k \in \Lambda }  2^{-\rho |k|}  . \label{1.31}
\eea
Let  $M \in \bbN$ be so large that $K \subset [-M+2, M-2] \times [-M+2, M-2]$. Then, $|n|^b \leq C_K|n +w|^b$ for all 
$w \in K $ and $n \in \Lambda$ with $|n| \geq 2 M$, hence,
\ef{1.24}, \ef{1.7}  and Remark \ref{rem1.2} imply
\bea
& & \sup_{w \in K, n\in \Lambda} \!\!\!\! \frac{|n|^{b} |\phi(w+n)|}{|\varphi(z)|}
\rowleq 
C_K \!\!\!\!    \sup_{w \in K , |n| < 2M} \!\!\!\!  \frac{|\phi(w+n)|}{|\varphi(z)|}
+  C_K \!\!\!\!  \sup_{w \in K, |n| \geq 2M } \!\!\!\!   |w+n|^{b} |\phi(w+n)| \leq C_K'.  \label{1.31b}
\eea
Since 
$$
\sum_{k \in \Lambda }   2^{-\rho |k|}  \leq \frac{C}{\rho^2},
$$
we get from \ef{1.31b} that \ef{1.31} is also bounded by $ C\rho^{-b}  (1 +  |m|)^{-b} .$
\ \ $\Box$

\bigskip

We will need in the next sections an invertible analytic multiplier function for switching the Floquet 
parameter in the quasiperiodic boundary conditions. Such a function can be constructed 
with the help of a second modification  of the elliptic function $\varphi$ in condition $(\cA)$.
The definition will be given in two steps. Recall that we  denote $m_1 = {\rm Re}\, m$, $m_2 = \ {\rm Im}\, m$ for all  $m \in \Lambda$.

\BEL  \label{lem1.4}
Let $\varphi$ and $\phi$ be as in \ef{1.6} and $\eta \in \bfQ$ and define  the functions 
\bea
\varphi_\eta (z) = \sum_{m \in \Lambda}  e^{-i \eta_1 m_1 - i\eta_2 m_2} \phi(z - m), \ z \in \Omega,   \ \ \ 
\mbox{and} \ \ \ {\widetilde \psi}_\eta  = \frac{\varphi_\eta}{\varphi}.   \label{1.32}
\eea

\smallskip \noindent
$(i)$ For every $\eta \in \bfQ$, the functions $\varphi_\eta$ and ${\widetilde \psi}_\eta$ 
are meromorphic in $\bbC$, and ${\widetilde \psi}_\eta$ is analytic on $\Omega$. 

\smallskip \noindent
$(ii)$
There exist constants  $\beta > 0$ and $C > 0$ such that 
%% function ${\widetilde \psi}_\eta$ is bounded and continuous on $\Omega$ and has, for some, the property
\bea 
\sup_{z \in \Omega, \eta \in \bfQ} | {\widetilde \psi}_\eta (z) | \leq C \ \ \ \ \mbox{and} \ \ \ \
\sup_{z \in \varpi} | {\widetilde \psi}_\eta (z)  - 1 | \leq C |\eta|^\beta \ \forall \, \eta \in \bfQ.   \label{1.32b}
\eea

\smallskip \noindent
$(iii)$ The function ${\widetilde \psi}_\eta$ satisfies the quasiperiodic boundary conditions 
\bea
{\widetilde \psi}_\eta(m + 1 + iy) = e^{i \eta_1} {\widetilde \psi}_\eta(m+iy), \ \ \ {\widetilde \psi}_\eta(m+x + i ) = e^{i \eta_2} {\widetilde \psi}_\eta(m+x)  \label{1.32a}
\eea
for all $\eta \in \bfQ$, $m \in \Lambda$ and $x,y \in (0,1) $. 
\ENL  

%%THE EXPONENT $\beta$ IS TO BE CONSIDERED. 

Proof. Assertion $(i)$ can be proved in the same way as Lemma \ref{lem1.3}.$(i)$, since the moduli of
the factors $e^{-i\eta_1 m_1 - i \eta_2 m_2}$ equal one.  As for $(ii)$, the first bound in \ef{1.32b} can be proved
in the same way as \ef{1.12h}. 

Let us prove the second inequality in \ef{1.32b}. We have for all $z \in {\rm cl} ( \varpi )  \cap \Omega$ and nonzero $\eta \in \bfQ$, by Remark \ref{rem1.2},
\beas
& &  | {\widetilde \psi}_\eta (z)  - 1 | \leq  \frac{\sum_{m\in \Lambda}  |1 - e^{-i \eta_1 m_1 - i\eta_2 m_2}| \, |\phi(z-m)| }{
\Big| \sum_{m\in \Lambda}   \phi(z-m) \Big| }
\rowleq 
 \frac{\sum_{|m| \leq  |\eta|^{-1/2} }  |1 - e^{-i \eta_1 m_1 - i\eta_2 m_2}| \, |\phi(z-m)| }{
\Big| \sum_{m\in \Lambda}   \phi(z-m) \Big| }
%C' \!\!\!\!\!
%\sum_{|m| \leq } \!\!\!\!\!  |1 - e^{-i \eta_1 m_1 - i\eta_2 m_2}| \, |\phi(z-m)| 
\rowpl 
C' \!\!\!\!\!
 \sum_{|m| > |\eta|^{-1/2} } \!\!\!\!\! |1 - e^{-i \eta_1 m_1 - i\eta_2 m_2}| \, |\phi(z-m)| =: S_1 + S_2 .
\eeas
To evaluate $S_1$ we note that, by the Taylor series of the exponential function,
\beas
|1 -  e^{-i \eta_1 m_1 - i\eta_2 m_2} | \leq C |\eta| \, |m| \leq C |\eta|^{1/2},
\eeas
which yields
\bea
S_1 \leq C \frac{|\eta|^{1/2} \sum_{m \in \Lambda  }  |\phi(z-m)| }{\big| \sum_{m \in \Lambda  }  \phi(z-m)\big| }
\leq C' |\eta|^{1/2}, 
\label{1.31t}
\eea
where we at the end again used an argument similar to that following \ef{1.28u}.

For $S_2$, we use $ |1 - e^{-i \eta_1 m_1 - i\eta_2 m_2}| \leq 2$ and an argument similar to \ef{1.11a}  to get 
\beas
& & S_2 \leq \!\!\!\!\!
\sum_{|m| > |\eta|^{-1/2} } \!\!\!\!\!  \, |m|^{2+\delta -b} |m|^{b- 2-\delta} |\phi(z-m)|
\rowleq 
|\eta|^{(b- 2-\delta)/ 2} \sum_{|m| > |\eta|^{-1/2}  }|m|^{b- 2-\delta} |\phi(z-m)| \leq  |\eta|^{(b- 2-\delta)/ 2} ;
\eeas
here, it is necessary to assume that $|\eta| < c_0$ for some small enough constant $c_0 >0$ so that, 
say, $|z- m| \geq 2$ holds for $z\in\varpi$ and $m$ with $|m| > |\eta|^{-1/2} > c_0^{-1/2}$. This is not a restriction, since we have already proved the 
first inequality in \ef{1.32b} and thus  the second one automatically holds for $|\eta| \geq c_0$, once
its constant $C$ is large enough. 

$(iii)$ The defining formula \ef{1.32} implies that $\varphi_\eta$ satisfies the quasiperiodic boundary conditions \ef{1.32a},
hence, also $\psi_\eta$ satisfies them, since $1/\varphi$ is doubly periodic.  \ \ $\Box$

\bigskip

As mentioned above, we will need an invertible multiplier operator changing the Floquet parameter values. Now, the second formula \ef{1.32b} implies
that $\widetilde \psi_\eta$ is an invertible function, if $|\eta|$ is small enough. However, there does
not seem be any guarantee that this is true for all $\eta \in \bfQ$: the function $\varphi_\eta$
and thus $\widetilde \psi_\eta$ could have zeros in the domain $\varpi$, which would destroy the 
invertibility of the associated multiplier operator. 

Let $R > 0$ be such that 
\bea 
\sup_{z \in \varpi} | {\widetilde \psi}_\eta (z)  - 1 | \leq \frac34 
 \ \ \ \mbox{for all} \ |\eta | \leq R ;   \label{1.33}
\eea
such an $R$ can always be found, by \ef{1.32b}. Clearly, \ef{1.33} implies
\bea
\inf_{z \in \varpi} |\widetilde \psi_\eta(z)| \geq \frac14 \ \ \ \mbox{for all} \ |\eta | \leq R. 
\label{1.34}
\eea
 In the case $R \geq \sqrt{2} \pi$, the inequality \ef{1.34}
holds for all $\eta \in \bfQ$, and we define $\psi_\eta = \widetilde \psi_\eta$ for all $\eta \in \bfQ$.
If $R < \sqrt{2} \pi$, we define
\begin{equation}
\psi_\eta (z) = \left\{
\begin{array}{ll}
\widetilde \psi_\eta(z) , \ \ & \mbox{for $\eta\in \bfQ$ with }|\eta| \leq R \\
\widetilde \psi_{R \eta /|\eta|} (z)^{|\eta|/R}, \ \ & \mbox{for $\eta\in \bfQ$ with }
  |\eta| > R. 
\end{array}
\right.  \label{1.36}
\end{equation}
If $z \in \varpi$, then $\widetilde \psi_{R \eta /|\eta|} (z) = r e^{i\theta}$ with $\theta
\in (-\pi/2, \pi/2)$, and the power in \ef{1.36} is defined as
$$
\widetilde \psi_{R \eta /|\eta|} (z)^{|\eta|/R} = r^{|\eta|/R} e^{i\theta|\eta|/R} .
$$ 
%%Obviously, the function  $\psi_\eta$ is analytic.

\BEL  \label{lem1.6}
Let $\eta \in \bfQ$. 

\smallskip \noindent
$(i)$ The function  $\psi_\eta$ is analytic  on $\Omega$, and there exists a constant $C > 0$ such that 
\bea 
\sup_{z \in \Omega, \eta \in \bfQ} |\psi_\eta (z)| \leq C .   \label{1.38k}
\eea

\smallskip \noindent
$(ii)$
If $\beta > 0$ is as in \ef{1.32b}, there holds %%, for some $\beta > 0$, the property
\bea 
\sup_{z \in \varpi} | 1 - \psi_\eta (z)| \leq \min\big( 3/4, C |\eta|^\beta \big)
\ \ \ \mbox{and} \ \ \ \inf_{z \in \varpi, \eta \in \bfQ} |\psi_\eta (z)| \geq 1/4  .   \label{1.38}
\eea

\smallskip \noindent
$(iii)$ The function $\psi_\eta$ satisfies the quasiperiodic boundary conditions 
\bea
\psi_\eta(m + 1 + iy) = e^{i \eta_1} \psi_\eta(m+iy), \ \ \  \psi_\eta(m+x + i ) = e^{i \eta_2}
\psi_\eta(m+x)  \label{1.40}
\eea
for all $m \in \Lambda$, $x,y \in (0,1)$. 
\ENL  

Proof. Given $\eta$, $\psi_\eta$ is analytic in $\Omega$ due to its definition \ef{1.36}, 
since all positive real powers of $\widetilde \psi_\eta$ can be defined as analytic functions, see the second inequality in \ef{1.32b}. 
All other  claims  in $(i)$ and $(ii)$  follow from Lemma \ref{lem1.4}, \ef{1.33}, \ef{1.34} and the definitions after them. 
Note that the upper bound $C |\eta|^\beta$ in \ef{1.38} only concerns $\eta$ in some neighborhood of 0, 
where $\psi_\eta$ and $\widetilde \psi_\eta$ coincide, hence \ef{1.38}  follows from \ef{1.32b}.

The quasiperiodicity \ef{1.40} follows from \ef{1.32a}, if $|\eta| \leq R$. If $|\eta| > R$,
we have
\beas
& & \psi_\eta(m+1 +iy) = {\widetilde \psi}_{R \eta /|\eta|} (m+1+iy)^{|\eta|/R} 
\roweq 
\big( e^{i R \eta_1 /|\eta|} \big)^{|\eta|/R}  {\widetilde \psi}_{R \eta /|\eta|} (m+iy)^{|\eta|/R} 
 = e^{i\eta_1} \psi_\eta(m+iy),
\eeas
for all $m \in \Lambda$ and $y \in (0,1) $ and, similarly, $\psi_\eta(m+ x +i)= e^{i \eta_2}   
\psi_\eta(m+x)y$ for all $x \in (0,1) $. \ \ $\Box$

\section{Mapping properties of the Floquet transform in Bergman spaces.}
\label{sec3}

We proceed to the first main result of the paper, namely, the characterization of  the image $\sF(A^2(\Omega))
\subset L^2(\bfQ; L^2(\varpi))$ in Theorem \ref{th2.3}. Naturally, this result is of importance, since 
one wants to preserve the bijectivity properties of $\sF$ also in the 
Bergman space case. Based on Theorem 3.6  in \cite{T1},  one cannot expect that just replacing $L^2(\varpi)$ by $A^2(\varpi)$ in the 
Bochner space makes the restriction of $\sF$ to $A(\Omega)$ a surjection: the quasiperiodic boundary conditions 
should appear in $\sF(A^2(\Omega))$. In this section we show  that with the help of the results in Section \ref{sec2}, the 
approach of \cite{T1} also works in the doubly periodic case. We will omit some details and refer 
to \cite{T1} for them. %% while still trying to keep the presentation of the approach self-contained. 

In fact, the range of the $\sF(A^2(\Omega))$ appears in the next definition. It involves the subspace
$A_{\eta, {\rm ext}}^2(\varpi)$ consisting of Bergman functions satisfying quasiperiodic boundary
conditions. Formula \ef{1.45} will give a lot of examples of such functions.

\BED
\label{def2.2}
If $\eta \in \bfQ$, we denote by $A_{\eta, {\rm ext}}^2(\varpi)$ the subspace of $A^2(\varpi)$ 
consisting of functions $f$ which can be extended as analytic functions in a neighborhood 
of ${\rm cl} (\varpi) \cap \Omega$ in $\Omega$ and satisfy the quasiperiodic 
boundary conditions 
\bea
f(1+i y)= e^{i \eta_1} f(iy), \ f(x +i)=e^{i\eta_2} f(x) \ 
\ \ \mbox{for all} \ x,y \in(0,1).  \label{2.7}
\eea
Moreover, $A_\eta^2(\varpi)$ stands for the closure of $A_{\eta, {\rm ext}}^2(\varpi)$ in $A^2(\varpi)$
and $ L^2(\bfQ; A_\eta^2(\varpi))$ for the subspace of $ L^2(\bfQ; A^2(\varpi))$ consisting of functions $f$ 
such that the function  $z \mapsto f (z, \eta)$ belongs to $A_\eta^2(\varpi)$ for a.e. $\eta \in  \bfQ$. 
\END

We continue with the following observations. 

\BEP
\label{lem1.1}
Assume condition $(\cA)$ holds for $\Omega$.

\noindent
$(i)$ For every $f \in A^2(\Omega)$, the function  $z \mapsto \sF f(z, \eta)$ is  analytic on $\varpi$ 
for a.e. $\eta \in \bfQ$. 

%\smallskip 
\noindent
$(ii)$ The space $ A^2(\Omega)$ has a dense subspace, denoted by $X$, which consists of functions 
$f$ such that, for all $\eta \in \bfQ$, the function $\sF  f ( \cdot, \eta)$ belongs to $A_ {\eta, {\rm ext}}^2(\varpi)$.
%%is analytic in 
%%$\varpi$ and has an analytic extension  to a neighborhood cl\,$(\varpi) \subset \Omega$.
\ENP

Proof.$(i)$ Let $\rho \in (0,1]$ and let $\psi^{(\rho)}$ be as in Lemma \ref{lem1.3}. We first show that for every 
$f \in A^2(\Omega)$ and $\eta \in \bfQ$, the function
\bea
\sF  ( \psi^{(\rho)} f)  (z ,\eta ) = 
\frac{1}{2\pi} \sum_{m \in  \Lambda} e^{- i \eta m } \psi^{(\rho)} (z + m) f(z+m)
\label{1.45}
\eea
is analytic in a neighborhood of $\varpi$. To this end, we define the set $U \subset \Omega$ to be the interior of 
the closure in $\Omega$  of the set
\beas
{\textstyle \bigcup}_{|m| \leq 2}  \varpi_m  
\eeas 
(note that the set $\varpi = \varpi_0$ is included in $U$). Then, let
$\xi \in U $ be arbitrary and pick up a small enough $\varrho > 0$ such that cl\,$ %\bea  
\big( B(\xi, \varrho) \big) \subset U. $ %\eea 
Applying a translation, the  Cauchy integral formula implies the estimate
\beas
\sup\limits_{z \in B(\xi, \varrho) }  |f%^{(n)} 
(z + m )| = \sup\limits_{z \in B(\xi, \varrho) +m} 
|f(z)| &\leq & \frac{C}{\varrho} \Vert f \Vert_{\varpi +m}
\leq 
\frac{C}{\varrho} \Vert f \Vert_{\Omega}   %\label{1.40}
\eeas
for all $m \in \Lambda$ and $z \in B(\xi , \varrho)$. The estimate \ef{1.15} implies that, for a fixed 
$\rho$, the series 
\ef{1.45} converges uniformly in the disc $B(\xi, \varrho)$. This means, 
$F  ( \psi^{(\rho)} f)$ is analytic in $B(\xi, \varrho)$ and thus in $U \supset \varpi$.

The rest of the proof is similar to that of Proposition 3.2 of \cite{T1}, except that we use
\ef{1.14} instead of (3.8)  of the reference. Indeed, \ef{1.14} and the unitarity of $\sF$  imply that 
$\sF  ( \psi^{(\rho)} f)$ converges to $\sF   f$ in $L^2(\bfQ;L^2(\varpi))$ as $\rho \to 0$,
and we thus find  a decreasing sequence  $(\rho_k)_{k=1}^\infty $ with $0 < \rho_k \to 0$ such that
\beas
\lim_{k \to \infty}  \sF  ( \psi^{(\rho_k)} f) (\cdot, \eta)  
= 
F   f (\cdot,\eta)  \ \ \mbox{in} \  L^2(\varpi)  %%\label{1.48}
\eeas
for almost all $\eta \in \bfQ$. %Now, the convergence in \eqref{1.48} happens in $A^2(\varpi)$ 
%for almost all $\eta \in \bfQ$, which 
This  implies that, for these $\eta$, the function
$\sF  f$ is analytic in $\varpi$, since $\sF  ( \psi^{(\rho)} f) (\cdot, \eta)$ is analytic in 
$\varpi$. 	

$(ii)$ The dense subspace $X$ can be defined to consist of all functions $\psi^{(\rho)} f$, where 
$f \in A^2(\Omega)$. As it was shown above, every  $\sF(\psi^{(\rho)} f)$ is an analytic function
in $U$, which contains a neighborhood of cl\,$(\varpi)$. It follows from the defining
formula \ef{1.45} that every 
$\sF(\psi^{(\rho)} f)$ satisfies the quasiperiodic conditions \ef{2.7}. Finally, the density of $X$
in $L^2(\Omega)$ follows from Lemma \ref{lem1.3}.$(iii)$.  \ \ $\Box$

\bigskip

We still consider some remarks and definitions, which are needed in proof of Theorem \ref{th2.3}.
First, we observe that the subspace   $ L^2(\bfQ; A_\eta^2(\varpi))$ of $ L^2(\bfQ; L^2(\varpi))$ is closed,   
since  $A_\eta^2(\varpi)$ is  closed  in $A^2(\varpi)$ for every $\eta$. We also define
$$
\cH_\eta := L^2 \big( \bfQ;  A_{\eta,{\rm ext}}^2(\varpi) \big)
$$ 
to be the subspace of $L^2(\bfQ; A_\eta^2(\varpi))$ which consists of functions $g$ 
such that the mapping $z \mapsto g(z ,\eta)$ belongs to 
$A_{\eta, {\rm ext}}^2(\varpi) $ for a.e. $\eta$. 

The following lemma can be proved in the same way as Lemma 3.5 of \cite{T1}, except that one uses the function $\psi_\eta$ of
Lemma \ref{lem1.6}.$(iii)$ instead of $e^{i \eta z}$ of the reference. Note that the proof is not completely trivial, since
the spaces   $ L^2(\bfQ; A_\eta^2(\varpi))$ and $\cH_\eta$ only  have the 
structure of a Banach vector bundle (e.g.  Section 1.3 of \cite{Kbook} or \cite{T1}, p. 209) and simple functions of 
$L^2 \big( \bfQ;  A^2(\varpi) \big)$ are not in general contained in these subspaces.

\BEL
\label{lem9}
The space $\cH_\eta$  is a dense subspace of $ L^2 \big( \bfQ;  A_{\eta}^2(\varpi) \big)$. 
\ENL

With these preparations we can state and prove the main result of this section. 

\BET
\label{th2.3}
Let $\Omega$ be doubly periodic domain satisfying assumption $(\cA)$.
The Floquet transform $\sF$ is a unitary operator from $A^2(\Omega)$ onto $ L^2(\bfQ; A_\eta^2(\varpi))$
with inverse  $\sF^{-1} :  L^2(\bfQ; A_\eta^2(\varpi)) \to A^2(\Omega)$ 
given by the formula \eqref{2.4}. 
\ENT

Proof. We first observe that by Proposition \ref{lem1.1}.$(ii)$,  $\sF$ maps the dense subspace $X$ of $A^2(\Omega)$ into $\cH_\eta$. 
In view of the unitarity of $\sF$, this implies that $\sF(A^2(\Omega))$ is contained in the closure of $\cH_\eta$
in $ L^2(\bfQ; L^2(\varpi))$, which equals $ L^2(\bfQ; A_\eta^2(\varpi))$, by Lemma \ref{lem9}.

We are thus left with the construction of an inverse image under $\sF $ for all $g \in   L^2(\bfQ; A_\eta^2(\varpi))$. 
Due to the density of $\cH_\eta$ in $L^2(\bfQ; A_\eta^2(\varpi))$ and the unitarity of $\sF$, we may assume that 
$g \in \cH_\eta$.
In the next construction we consider $\eta \in \bfQ$ such that the function
$g(\cdot , \eta)$ is a well-defined element of $A_{\eta,{\rm ext} }^2(\varpi)$. %%[[BETTER]]
For all  such $\eta \in \bfQ$, we define the function $G_\eta:\Omega \to \bbC$ by setting
\bea
G_\eta(z) = G_{\eta,m} (z) \ \ \ \mbox{for all} \ z \in {\rm cl} ( \varpi_m), \ m \in \Lambda, \label{2.25}
\eea
where 
\bea
G_{\eta,m} : \varpi_m \to \bbC, \ \ G_{\eta,m}(z) =  e^{i \eta_1 m_1 + i \eta_2 m_2} 
g( z -  m , \eta).   \label{1.60}
\eea
Let us show that $G_\eta$ is well defined and gives an analytic function in $\Omega$. 
Indeed, $G_{\eta,m}$ 
is analytic on the subdomain $\varpi_m$, $m\in \Lambda$ (see \eqref{1.1})
and  has an analytic extension, still denoted by $G_{\eta,m}$, to a  neighborhood  
of the closure of $\varpi_m$ in $\Omega$ (see Definition \ref{def2.2}). 
Moreover,  the functions $g(\cdot ,\eta)$ %%and $\psi_\eta$ 
satisfy the quasiperiodic boundary conditions  \eqref{2.7}. %%, see Lemma \ref{lem1.4} for the latter. 
Hence, we get for all 
$m \in \Lambda$, all $z = m +1 + iy$ with  $y \in (0,1) $,
\beas
& &G_{\eta , m} (z) =  e^{i \eta_2 m_2} e^{i \eta_1 m_1}  g(z- m,\eta ) =   e^{i \eta_2 m_2} e^{i \eta_1 m_1}  g(1 +iy,\eta )
\roweq  e^{i \eta_2 m_2}e^{i (m_1 + 1)  \eta_1  } g(iy , \eta ) 
 =  e^{i \eta_2 m_2} e^{i (m_1 + 1)  \eta_1  } g(z -(m+1) , \eta ) %%=  G_{\eta,m+1} (m +1  +iy )
 = G_{\eta, m+1} (z).
\eeas
%for all  $y \in [a,b]$. 
We see that  the functions $G_{\eta,m}$ and $G_{\eta,m+1}$ 
%%are two analytic functions with overlapping domains of definition  and they  
coincide on  the line  segments of positive length, which consist of the common parts of
$\partial \varpi_m$ and $\partial \varpi_{m+1}$. Due to the their analyticity,  $G_{\eta,m}$ and 
$G_{\eta,m+1}$ thus coincide in their entire common domain of definition.
In the same way one shows that the functions $G_{\eta,m}$ and $G_{\eta, m+i}$ coincide in their (non-empty) 
common domain of definition. Since these claims hold for all $m \in \Lambda$, we conclude that the 
function $G_\eta$ in \ef{2.25} is well-defined and analytic in $\Omega$. 

From  \eqref{1.60}, \eqref{2.4} we observe that 
\beas
H:  =  \frac{1}{2\pi}\int\limits_{\bfQ} G_\eta (z) d\eta  =  \sF^{-1} g (z)\in A^2(\Omega),
\eeas
thus, there holds $\sF H = g$ for all $g \in \cH_\eta \subset L^2\big(\bfQ;
L^2(\varpi)\big)$. %%, by the general Floquet inversion formula \eqref{1.4}. 
We conclude that image of $\sF \big( A^2(\Omega) \big) $ contains the subspace $\cH_\eta$. In 
view of the fact that ${\sf F}$ is an isometry,  $\sF \big( A^2(\Omega) \big) $ also contains 
the entire space %%there remains to show that $\cH_\eta$ is dense in 
$L^2(\bfQ; A_\eta^2(\varpi))$. This completes the proof. \ \ $\Box$

\section{General kernel formula for the periodic domain.}
\label{sec4}

The results of the previous section identify the $\eta$-dependent family of Bergman spaces $A_\eta^2 (\varpi)$
on the periodic cell $\varpi$, which corresponds, under the Floquet transform, to the Bergman space on the 
original doubly periodic domain $\Omega$. The quasiperiodic boundary conditions in the periodic cell do 
not show up in the corresponding $L^2$-isometry (see the beginning of Section \ref{sec1}), but the 
situation in the Bergman space case is analogous to the Sobolev space case, see \cite{Ku}, Section 4.2. 
Let us next present the connection of the Bergman-type projections in the domains  $\Omega$ and $\varpi$ 
by using the Floquet transform. The results and their proofs are completely analogous to the singly
periodic case so that a repetition of the proofs will not be necessary.

If  $\eta \in  \bfQ$, let us  denote  the orthogonal projection from $L^2 (\varpi)$ onto 
$A_{\eta }^2 (\varpi)$ by $P_\eta$. The projection $P_\eta$ can be written as an integral
operator 
\bea
P_\eta f (z) = \int\limits_\varpi K_\eta(z,w) 
f(w) dA(w) , \label{1.17}
\eea
since the existence of kernel can be proved following the usual proof for the existence of 
Bergman kernels, see  \cite{K1},  p.1060. In particular, the kernel has the property 
that $K_\eta(z, \cdot) \in L^2(\varpi)$ for all $z \in \varpi$. 

We next  define the bounded operator $\cP:  L^2(\bfQ; L^2(\varpi)) \to L^2(\bfQ; A_\eta^2(\varpi))$ by 
\bea
\cP f(z,\eta) = \big( P_\eta f( \cdot , \eta ) \big)  (z) , \ \ \ 
f \in L^2(\bfQ; L^2(\varpi)), \ z \in \varpi, \ \eta \in \bfQ.  
\eea
The following result can now be proved in the same way as  Lemma 4.2 and Theorem 4.3 in  \cite{T1}.
We denote here $z^\tr =  z -  [{\rm Re} z] -  i [{\rm Im} z]$ so that $z^\tr \in {\rm cl}( \varpi)$ 
for all $z \in \Omega$.

\BET
\label{th3.3} 
Let $\Omega$ be doubly periodic domain satisfying assumption $(\cA)$.

\noindent $(i)$ The operator $\cP$ is the orthogonal projection from
$L^2(\bfQ; L^2(\varpi))$  onto $L^2(\bfQ; A_\eta^2(\varpi))$. 

\noindent $(ii)$ 
The Bergman projection $P_\Omega$ from $L^2(\Omega)$ onto $A^2(\Omega)$ can be written as  
\bea
& &  \sF^{-1} \cP  \sF  f (z) = 
\frac1{2 \pi} \int\limits_\bfQ e^{ i [{\rm Re} z ]\eta_1 + i [{\rm Im} z ]\eta_2 }
\big( P_\eta  \widehat f ( \cdot , \eta) \big) (z^\tr  )  
d \eta 
\roweq
\frac1{4 \pi^2} \int\limits_{\Omega} \int\limits_\bfQ 
e^{ i \eta_1 ( [{\rm Re} z]-  [{\rm Re} w ]) +  i \eta_2 ( [{\rm Im} z ] - [{\rm Im} w ] ) }  %\times
%\nonumber \\
%& & \hskip1.5cm \times 
K_\eta( z^\tr  , w^\tr) 
f(w) d\eta d A(w)  
\label{4.4zz}
\eea
\ENT

\section{Spectral band formula for Toeplitz operators with doubly periodic symbols}
\label{sec5}

Given a doubly periodic domain $\Omega$ and a symbol $a \in L^\infty(\Omega)$, the Toeplitz operator 
$T_a : A^2(\Omega) \to A^2(\Omega) $ is defined by 
$$
T_a f=P_\Omega M_a f=P_\Omega(af) 
$$
where $M_a$ is the pointwise multiplier $f \mapsto af$, $f \in A^2(\Omega) $. The aim of this section is
to study the spectra $\sigma(T_a)$ and essential spectra $\sigma_{\rm ess} (T_a)$ of the operators 
$T_a$ with doubly periodic symbols $a\in L^\infty(\Omega)$, 
\begin{equation}\label{periodic}
a(z)=a(z+1)=a(z+i) \text{ for almost all } z\in\Omega . 
\end{equation}
In the main result, Theorem \ref{th8.3}, we show that the spectrum and essential spectrum of such a $T_a$ 
coincide and, moreover, can be presented as the union of the spectra of a family of 
Toeplitz-type operators $T_{a,\eta}:A_{\eta}^2(\omega)\to A_{\eta}^2(\omega)$,
$\eta \in \bfQ$, in the periodic cell. These operators are defined with the help of the projections $P_\eta$ of \ef{1.17} by 
$$
T_{a,\eta}f=P_{\eta}(a|_{\varpi}f)
$$
for all $\eta \in \bfQ$, $f \in A_{\eta}^2(\omega)$, and their spectra as operators in $A_{\eta}^2(\omega)$
are denoted by $\sigma(T_{a,\eta})$.

The main result reads as follows. 

\BET  \label{th8.3}
Let $\Omega$ be doubly periodic domain with property $(\cA)$ and let the symbol $a \in L^\infty(\Omega)$
be as in \eqref{periodic}.
The essential spectrum of the Toeplitz-operator $T_a:A^2(\Omega) \to A^2(\Omega)$ can be described by the
formula
\bea
\sigma_{\rm ess} ( T_a )  = \bigcup_{\eta \in \bfQ} \sigma(T_{a,\eta}). 
\label{8.40}
\eea
Moreover, there holds $\sigma( T_a )= \sigma_{\rm ess} ( T_a )$.
\ENT

The proof follows that of Theorem 4.1. of \cite{T2}, but the main difference is the need to use
the  multiplier $J_{\mu,\eta}$, the definition of which is based on the much more delicate 
function $\psi_\eta$ of Lemma \ref{lem1.6} instead of the basic exponential function $e^{i\eta z}$ of 
\cite{T2}. This causes some differences to the arguments. For example, the function $1/\psi_\eta$ is not the 
same as $\psi_{-\eta}$ as in the case of the exponential function, which could complicate
the considerations involving the inverse operators, and the  estimate \ef{4.4} in below is
worse than its analogue (3.9) in \cite{T2}. It is thus worthwhile to present the main steps of the proof. 

We extend the definition of the operators $T_{a,\eta}$ to the Bochner space by defining  $\mathcal{T}_a: L^2(\bfQ;A_{\eta}^2(\varpi))\to L^2(\bfQ;A_{\eta}^2(\varpi)) $ as 
$$
\mathcal{T}_a :f(\cdot,\eta)\mapsto  T_{a,\eta}f(\cdot,\eta). 
$$
Also, we denote by $\mathcal{M}_a:L^2(\bfQ;A_{\eta}^2(\varpi))\to L^2(\bfQ;A_{\eta}^2(\varpi))$ 
the operator 
$$
\mathcal{M}_a :f(\cdot,\eta)\mapsto  a|_{\varpi}f(\cdot,\eta)
$$
so that there holds $\mathcal{T}_a =\mathcal{P}\mathcal{M}_a $. The definitions easily imply 
$
T_af=\sF^{-1}\mathcal{T}_a\sF f
$  for all $f\in A^2(\Omega)$.

\bigskip

We assume that $(\cA)$ holds for $\Omega$ and refer to Lemma \ref{lem1.6} for  the definition of the 
analytic function  $\psi_\eta$. Given $\eta, \mu \in\bfQ$, we define the operators
\bea
J_{\eta, \mu} f(z)= \psi_{\mu-\eta} (z)  f(z) \ \ \ \mbox{and} \ \ \ 
T_{a, \eta, \mu}=J_{\eta, \mu}^{-1} T_{a, \mu} J_{\eta, \mu}. %%=J_{\mu, \eta} T_{a, \mu} J_{\eta, \mu}
\label{4.2}
\eea

\BEL  \label{lem12}
Given $\eta, \mu \in \bfQ$, the operator  $J_{\eta, \mu}$ is a bounded bijection of $A^2(\varpi)$ onto  itself and 
also a bounded bijection from $A_\eta^2(\varpi)$ onto $A_\mu^2(\varpi)$.
Consequently, $T_{a, \eta, \mu}$ is a bounded operator from $ A_\eta^2(\varpi)$ into itself
and from $A^2(\varpi)$ into itself. 
\ENL 
 
Proof. 
It follows from Lemma \ref{lem1.6}.$(i),(ii)$ that $J_{\eta, \mu}$ is a bounded bijection of $A^2(\varpi)$ onto  itself, and its 
inverse is, of course, the multiplication operator $J_{\eta, \mu}^{-1} : f \mapsto f/ \psi_{\mu -\eta}$.
Also, by Lemma \ref{lem1.6}.$(iii)$, if $f \in A_{\eta, {\rm ext}}^2(\varpi)$, then the function $\psi_{\mu-\eta} f$
satisfies the quasiperiodic conditions \ef{2.7} with the Floquet parameter $\mu \in \bfQ$. Thus, 
$J_{\eta, \mu}$ maps  $A_\eta^2(\varpi)$ boundedly into $A_\mu^2(\varpi)$. Also, 
$J_{\eta, \mu}^{-1}$ maps  $A_\mu^2(\varpi)$  into $A_\eta^2(\varpi)$, since the function $1/\psi_{\mu -\eta}$ satisfies
the quasiperiodic conditions with the parameter $\eta- \mu$. We conclude that $J_{\eta,\mu}$ is 
also a bounded bijection from $A_\eta^2(\varpi)$ onto 
$A_\mu^2(\varpi)$. 
\ \ $\Box$

\bigskip

The operators $J_{\eta,\mu}$ are not isometries, since 
the factor $\psi_\eta$ is not unimodular. However, \ef{1.38} implies for all $f \in L^2(\varpi)$,
\beas
\int\limits_\varpi |f - \psi_\eta f|^2 dA = 
\int\limits_\varpi |1 - \psi_\eta |^2|f|^2  dA  \leq C|\eta|^{2 \beta} \Vert f\Vert_\varpi^2, 
\eeas
hence, 
\bea
\left\|I_{\varpi}-J_{\eta, \mu}\right\|_{L^2(\varpi) \rightarrow L^2(\varpi)} \leq C|\eta-\mu|^\beta 
\ \ \  \forall \mu, \eta \in\bfQ   \label{4.4}
\eea
where $I_{\varpi}$ is the identity operator on $L^2(\varpi)$. Also, Lemma \ref{lem1.6}.$(i),(ii)$ implies 
\beas
\sup_{z \in \varpi} \Big| 1- \frac{1}{\psi_\eta(z) }\Big| = \sup_{z \in \varpi}  \frac{| 1 - \psi_\eta(z) |}{|\psi_\eta(z)| }
\leq  \frac{ \sup_{z \in \varpi} | 1 - \psi_\eta(z) |}{ \inf_{z \in \varpi}|\psi_\eta(z)| } \leq C|\eta|^\beta
\eeas
and consequently, 
\bea
\left\|I_{\varpi}-J_{\eta, \mu}^{-1}\right\|_{L^2(\varpi) \rightarrow L^2(\varpi)} \leq C|\eta-\mu|^\beta 
\ \ \  \forall \mu, \eta \in\bfQ .   \label{4.4b}
\eea
We also get, by the definition of $T_{a, \eta, \mu}$, 
$$
\|T_{a,\mu}-T_{a, \eta, \mu}\|_{A^2(\varpi)\to A^2(\varpi)}\leq C |\mu-\eta|^{\beta}.
$$

\BEL  \label{lem8.0} 
There exist constants $C, C' >0$ such that if  $\eta, \mu \in \bfQ$, we have 
$\Vert P_\eta - P_\mu \Vert_{L^2(\varpi) \to L^2(\varpi)} \leq C |\eta- \mu|^{\beta/2}$ 
and consequently $\Vert T_{a,\eta} - T_{a,\mu} \Vert_{A^2(\varpi) \to A^2(\varpi)} \leq 
C' |\eta- \mu|^{\beta/2} $. 
\ENL

Proof.  If  $\mu \in \bfQ$ is fixed, the idea of the proof, according to \cite{T1}, 
is to consider the non-orthogonal projections $\widetilde P_\mu = J_{\mu,\eta}^{-1} P_\eta J_{\mu,\eta}$  
from $L^2(\varpi) $ onto $A_\mu^2(\varpi)$, where $\eta \in \bfQ$. 
For every  $f \in L^2(\varpi)$ we denote $f_A = P_\mu f$, $f^\perp =
f -f_A$, and observe that $(P_\mu - \wP) f = (P_\mu - \wP) f^\perp$, since 
both $P_\mu$ and $\wP$ project onto $A_\mu^2(\varpi)$. This yields 
\bea
& & \big( (P_\mu - \wP) f \big| (P_\mu - \wP) f \big)_\varpi 
=  \big( (P_\mu - \wP) f^\perp \big| (P_\mu - \wP) f^\perp \big)_\varpi
\roweq
( \wP f^\perp | \wP f^\perp )_\varpi 
 =  ( \wP f^\perp - f^\perp | \wP f^\perp )_\varpi ,
\label{8.s6}
\eea
where we at the end used $ (  f^\perp |\wP f^\perp )_\varpi=0$, which follows from  $ \wP  f^\perp \in 
A_\mu^2(\varpi)$ and $f^\perp \in A_\mu^2(\varpi)^\perp$. 
We will soon show that 
\bea
\big| ( \wP f^\perp - f^\perp | \wP f^\perp )_\varpi  \big| 
%%\Vert \wP f^\perp - f^\perp \Vert^2 
\leq C |\eta -\mu|^\beta \Vert f \Vert_\varpi^2.  \label{8.s12}
\eea
Combining \ef{8.s6}--\ef{8.s12} yields $\Vert P_\mu - \wP \Vert_{L^2(\varpi) \to L^2(\varpi)}$ $ \leq 
C |\eta -\mu|^{\beta/2}$, which proves the result, since we obtain from \ef{4.4}, \ef{4.4b}
\beas
& & \Vert P_\eta -P_\mu \Vert_{L^2(\varpi) \to L^2(\varpi)} \leq 
\Vert P_\eta -\wP \Vert_{L^2(\varpi) \to L^2(\varpi)} + \Vert \wP - P_\mu  \Vert_{L^2(\varpi) \to L^2(\varpi)}
\rowleq 
\Vert( I_\varpi -  J_{\mu,\eta}^{-1}) P_\eta \Vert_{L^2(\varpi) \to L^2(\varpi)}  
+ %% \rowpl
\Vert J_{\mu,\eta}^{-1} P_\eta (  I_\varpi - J_{\mu,\eta} ) \Vert_{L^2(\varpi) \to L^2(\varpi)}  
+C |\eta -\mu|^{\beta/2}
\rowleq C' |\eta -\mu|^{\beta/2}.
\eeas

To see \ef{8.s12}, \ef{4.4b} implies
\beas
\Vert P_\eta J_{\mu,\eta} f - \wP f \Vert_\varpi
\leq C |\eta - \mu|^\beta \Vert f \Vert_\varpi
%\label{8.s2}
\eeas
for all  $f \in L^2(\varpi)$. This and \ef{4.4} again yield for  $f \in L^2(\varpi)$
\beas
& & \big| \big( f - \wP f \big| \wP f \big)_\varpi \big|
\rowleq 
\Big| \big( f - J_{\mu,\eta}f +P_\eta J_{\mu,\eta} f- \wP f \big| \wP f \big)_\varpi \Big|
+%\rowpl
\Big|  \big( J_{\mu,\eta} f - P_\eta J_{\mu,\eta} f \big| \wP f - P_\eta J_{\mu,\eta} 
f \big)_\varpi  \Big| 
\rowleq
\Big| \big( f - J_{\mu,\eta}f  \big| \wP f \big)_\varpi \Big|
+ \Big| \big( P_\eta J_{\mu,\eta} f- \wP f \big| \wP f \big)_\varpi \Big|
+ %%\rowpl  
C |\eta - \mu|^\beta \Vert f \Vert_\varpi^2 
\rowleq
C' |\eta - \mu|^\beta \Vert f \Vert_\varpi^2 . \ \ \Box  %%\label{8.s4}
\eeas

\bigskip

Continuing the proof of Theorem \ref{th8.3}, the proof of Lemma 4.2 of \cite{T2} shows that 
the set 
$
\Sigma := \bigcup_{\eta \in \bfQ} \sigma(T_{a,\eta})  %\label{8.14a}
$ 
is closed. Namely, if the claim were not true, one could find  $\lambda \in 
{\rm cl}\,(\Sigma) \setminus \Sigma$ and sequences $(\eta_k)_{k=1}^\infty 
\subset \bfQ$ and $(\lambda_k)_{k=1}^\infty$ such that $\lambda_k \in 
\sigma (T_{a,\eta_k} )$ and $\lambda_k \to \lambda$ as $k \to \infty$. 
By passing to a subsequence, if necessary, one  may assume 
that  $\eta_k \to \eta$ for some $\eta \in \bfQ$ as $k \to  \infty$. 
Moreover, since $\sigma(T_{a,\eta})$ is closed, there is a number $0< \delta <1$ such that 
$ %\bea
{\rm dist}\,( \lambda, \sigma(T_{a,\eta}) ) \geq \delta \ (*) . 
$ 
Then, one observes that, due to Lemma \ref{lem8.0},
\bea
& & \Vert T_{a,\eta_k,\eta} - T_{a,\eta_k} \Vert_{A^2(\varpi) \to A^2(\varpi)} 
\leq C |\eta -\eta_k|^{1/2} .  \label{8.20q}
\eea
Moreover,  the spectra of $T_{a,\eta}$  in  $A_\eta^2(\varpi)$ and $ T_{a,\eta_k,\eta}$ in 
$A_{\eta_k}^2(\varpi)$ coincide (cf. Lemma 3.3. of \cite{T2}), 
which implies that, for a large enough $k$,  $\lambda_k$ belongs to the resolvent set of 
$T_{a,\eta_k,\eta}$, see $(*)$. This eventually leads to a contradiction with $\lambda_k 
\in \sigma (T_{a,\eta_k} )$ by using some quite standard
resolvent estimates and the closeness of the two operators given by formula \ef{8.20q}.

The closedness of $\Sigma$ implies that $\sigma(T_a) \subset \Sigma$. To see this,
we fix a number $\lambda \in \bbC$ with $\lambda \notin \Sigma$. 
Then, for every $\eta \in \bfQ$, there exists a bounded inverse
$R_{\eta, \lambda} : A_\eta^2(\varpi) \to A_\eta^2(\varpi)$ of the operator 
$T_{a,\eta} - \lambda I_\varpi$. Since  the operator norm 
$\Vert R_{\eta, \lambda} \Vert_{ A_\eta^2(\varpi) \to  A_\eta^2(\varpi)}$  depends continuously on $\eta$
(see the proof of Corollary 4.3 in \cite{T2}), it has  a uniform
upper bound for all $\eta \in \bfQ$. 
We conclude that the operator $\cR_\lambda : f(\cdot ,\eta) \mapsto R_{\eta, \lambda} 
f(\cdot,\eta)$ is bounded in $L^2(\bfQ ; A^2_\eta(\varpi))$, and it is the inverse of the operator
$\cT_a - \lambda \cI$ in the space  $L^2(\bfQ ; A^2_\eta(\varpi))$, where 
$\cI$ is the indentity operator on $L^2(\bfQ ; A^2_\eta(\varpi))$. Hence, 
$\sF^{-1} \cR_\lambda \sF$ is a bounded inverse of $T_a - \lambda 
I_\Omega = \sF^{-1} ( \cT_a - \lambda\cI) \sF$, and $\lambda$ thus belongs to the 
resolvent set of $T_a$. We refer to \cite{T2} for the details of these arguments. 

\bigskip

To complete the proof of Theorem \ref{th8.3}, it is thus sufficient to show that
$\Sigma \subset \sigma_{\rm ess} (T_a)$, i.e. every $\lambda$, which belongs to $\sigma(T_{a,\mu})$ for some 
$\mu \in \bfQ$, is a point in the essential spectrum of $T_a$. To this end, we will use the Weyl criterion
as presented e.g. in Lemma 1.1. of \cite{T2}, and thus we  fix $\mu \in \bfQ$, $\lambda \in 
\sigma(T_{a,\mu})$, an arbitrary  $\varepsilon > 0$ and a near eigenfunction $g \in A_{\mu}^2(\varpi)$ with 
$\Vert g \Vert_\varpi = 1$ such that
\bea
\Vert T_{a,\mu} g - \lambda g \Vert_\varpi \leq \varepsilon. \label{8.42}
\eea
Since the operator $T_{a,\mu}$ is bounded and $A_{\mu, {\rm ext}}^2(\varpi)$ is dense in 
$A_{\mu}^2(\varpi)$, it can be assumed that $g \in A_{\eta, {\rm ext}}^2(\varpi)$.

Given $n \in \bbN$, we define  
$ %\bea
G_n = \sF^{-1} (\cX_{\mu,n} g ) ,   %\label{8.43}
$ % \eea
where  $  (\cX_{\mu,n}g) (z,\eta) := g \otimes\cX_{\mu,n}  (z,\eta) =  
g(z)\cX_{\mu,n}(\eta) \in L^2(\bfQ;A_\eta^2(\varpi))$ with 
\begin{equation}
\cX_{\mu,n}(\eta) = \left\{
\begin{array}{ll}
n , \ \ \ & \eta \in B(\mu, 1/n) \\
0, &\mbox{for other} \ \eta \in \bfQ, 
\end{array}
\right. \label{8.43a}
\end{equation}
where $B(\mu, 1/n) = \{ \eta \in \bfQ \, : \, |\eta - \mu| < 1/n \}$.
We have $1/2 \leq \int_\bfQ \cX_{\mu,n} (\eta)^2 d \eta \leq 4$ for all $n$, hence 
\beas
1/2 \leq 
\Vert \cX_{\mu,n} g \Vert_{L^2 (\bfQ; L^2(\varpi)) }^2
= \int\limits_\bfQ \cX_{\mu,n} (\eta)^2 \Vert g \Vert_\varpi^2  d \eta \leq 4 .
\eeas
Taking the Floquet inverse transform yields 
$G_n = \sF^{-1} ( \cX_{\mu,n} g )  \in  A^2(\Omega)$ and 
\beas
1/2 \leq \Vert G_n  \Vert_\Omega \leq 4.  %%\label{8.43p}
\eeas
Using   Lemma \ref{lem1.3} we choose for every $n$ the number $\rho(n) > 0$ 
such that $\rho(n+1) < \rho(n)$ and 
%such that $\Vert G_n -   G_n  \Vert_\Omega \leq 1/(n+2)$. By 
%reindexing  $\varphi_{\ell(n)} \to \varphi_n$ we thus have, for all $n \in \bbN$,   
\bea
\Vert G_n -   \psi^n  G_n  \Vert_\Omega \leq \frac{1}{n +2}  \ \ \ \mbox{with}
\ \ \ \psi^n := \psi^{(\rho(n))}.
\label{8.43pp}
\eea

The Weyl singular sequence for the number $\lambda$ is defined by the translated functions $h_n = f_n \circ \sft_n
\Vert f_n\Vert_\Omega^{-1}$, where   $ % \bea
f_n = \psi^n G_n \in A^2(\Omega) 
$ % \eea
and $\sft_n(z) = z -  m(n) $, and the numbers $m(n)$ are chosen so as to form a fast enough
increasing sequence of positive integers. 
Note that the above definitions imply  $1/C' \leq \Vert f_n \Vert_\Omega \leq C'$ for some constant $C'>0$, for all $n$,
hence, it is plain that condition (1.4) of Lemma 1.1 in \cite{T2} holds, once we  verify the  estimate 
\bea
\Vert T_a f_n - \lambda f_n \Vert_\Omega \leq C \varepsilon \label{8.39k} 
\eea 
for %%any arbitrary $\varepsilon > 0$ and 
all large enough $n \in \bbN$ (here, the constant $C$ may depend on $\lambda$). To see this, %%for the given $\varepsilon$ 
we assume that $n  \geq \varepsilon^{-2/\beta}$ in the following. 
Since 
$ %\bea
\Vert T_a ( f_n -  G_n) \Vert_\Omega 
< C \varepsilon %%label{8.56a}
$ %\eea
by \ef{8.43pp} and the boundedness of the operator $T_a$ in $L^2(\Omega)$, it suffices to show that 
$ 
\Vert T_a  G_n - \lambda G_n \Vert_\Omega < C \varepsilon . 
$

We write 
\beas
T_a G_ n &= &  \sF^{-1} \cP \cM_a (\cX_{\mu,n} g  ) 
= \sF^{-1}   \big(  \cX_{\mu,n }(\eta) P_\eta ( a g  )(z)  \big) 
\roweq
\sF^{-1}  \Big(   \cX_{\mu,n } (\eta) \big( P_\eta ( a g  ) (z)
- P_\mu (ag)(z) \big) \Big)  
\rowpl
\sF^{-1}   \big( \cX_{\mu,n } (\eta) P_\mu (ag)(z)  - \lambda  \cX_{\mu,n }(\eta) g(z) \big)  
+ \sF^{-1}   \big( \lambda  \cX_{\mu,n } g \big)  
\nonumber \\
&=:& 
\Psi_1 + \Psi_2  +\lambda G_n .  %\label{8.56}
\eeas
Lemma \ref{lem8.0}  yields the bound   $\Vert P_\eta - P_\mu \Vert_{A^2(\varpi) \to 
A^2 (\varpi)} \leq C|\eta- \mu|^{\beta/2}$. We get by \ef{8.43a} %%of the symbol $a$, 
\bea
& & \big\Vert    \cX_{\mu,n } \big( P_\eta ( a g  ) - 
P_\mu (ag)  \big) \big\Vert_{L^2(\bfQ;L^2(\varpi))}^2
\rowleq 
C_a \!\!\!\!\! \int\limits_{B(\mu, 1/n)} \!\!\!\!\!\! 
n^2 |\eta -\mu|^\beta \Vert g \Vert_\varpi^2 d\eta
\leq %\frac{C'}{n} \Vert g \Vert_\varpi^2 = 
\frac{C'_a}{n^\beta} \leq C'_a \varepsilon^2 .   \label{8.56w}
\eea  
The near eigenfunction property \ef{8.42} yields a similar estimate for $\Psi_2$:
\beas
& & \big\Vert    \cX_{\mu,n } \big( P_\mu ( a g  ) - 
\lambda g  \big)\big\Vert_{L^2(\bfQ;L^2(\varpi))}^2
\leq 
C  \!\!\!\!\!\! \int\limits_{B(\mu,1/n)} \!\!\!\!\!\!
n^2 \varepsilon^2  d\eta \leq C' \varepsilon^2. 
\eeas
From this and \ef{8.56w} we obtain $\Vert \Psi_1 \Vert_\Omega + 
\Vert \Psi_2 \Vert_\Omega \leq C \varepsilon$, and this shows that the functions
$h_n$ satisfy \ef{8.39k} and thus  (1.4) of Lemma 1.1 in \cite{T2}.

There remains to show that the sequence $(h_n)_{n=1}^\infty$ does not have subsequences converging
in $L^2(\Omega)$. 
This is intuitively quite clear, since the functions $f_n$ have some localization in $\Omega$ 
due to the factors $\psi^n$. Choosing a rapidly enough growing sequence $(m(n))_{n=1}^\infty$
for the definition of the translations $\sft_n$, see above, one can arrange that
$\Vert h_n - h_m \Vert_\omega$ has a positive lower bound independent of $n,m$, since
$h_n$ is essentially the same as $f_n \circ \sft_n $. The details of the proof can be completed 
by replacing $\varphi_n$ by $\psi^n$ and by other obvious changes at the end of Section 4 of \cite{T2}.  
This completes the proof of Theorem \ref{th8.3}.

\end{document}